\documentclass{amsart}
\newtheorem{theorem}{Theorem}

\newtheorem{remark}[theorem]{Remark}

\newtheorem{corollary}[theorem]{Corollary}

\newtheorem{definition}[theorem]{Definition}
\newtheorem{example}[theorem]{Example}

\newtheorem{lemma}[theorem]{Lemma}

\renewcommand{\oddsidemargin}{0.5cm}

\begin{document}
\title[Perron Theorem and Monotone Iteration Method]{Perron Theorem in the Monotone Iteration Method for Traveling Waves in Delayed
Reaction-Diffusion Equations}
\author{Amin Boumenir}
\address{Department of Mathematics, University of West Georgia, Carrollton, GA 30118}
\email{boumenir@westga.edu}
\author{Nguyen Van Minh}
\address{Department of Mathematics, University of West Georgia, Carrollton, GA 30118}
\email{vnguyen@westga.edu}
\thanks{The authors thank the referee for reading carefully the
manuscript, and for suggestions to improve the presentation of the
paper.}
\date{\today}
\begin{abstract}In this paper we revisit the existence of traveling waves for delayed reaction
diffusion equations by the monotone iteration method. We show that Perron
Theorem on existence of bounded solution provides a rigorous and constructive
framework to find traveling wave solutions of reaction diffusion systems with
time delay. The method is tried out on two classical examples with delay: the
predator-prey and Belousov-Zhabotinskii models.
\end{abstract}
\keywords{Traveling wave, Perron Theorem, monotone iteration, predator prey model,
Belousov-Zhabotinskii model}
\subjclass{35K55, 35R10}
\maketitle

\section{Introduction}

We shall be concerned with the existence of traveling waves solutions for the
delayed reaction diffusion system
\begin{equation}
\frac{\partial u(x,t)}{\partial t}=D\frac{\partial^{2}u(x,t)}{\partial x^{2}%
}+f(u_{t}).\label{pde-0}%
\end{equation}
Due to their important applications in population dynamics and biological
models, see e.g.
\cite{mur,gar,gar2,3vol,dun1,dun2,mishut,misrei,wu,sch,zouwu,wuzou,ma,huazou,hualurua,wanlirua,zha}%
, equations such as (\ref{pde-0}) have evolved from the simple one
dimensional scalar reaction diffusion equation to systems that
include a delay in time for more realistic modeling. Traveling
waves, although are a classical topic in applied mathematics,
remain a driving force in the study of (\ref{pde-0}). We refer the
reader to \cite{gar3,3vol} which contain many surveys on methods
used to study traveling waves in parabolic differential equations.
Note that when delay is introduced, most methods would fail if
they are not modified appropriately.

\medskip Recently, Wu and Zou \cite{zouwu,wuzou} have extended the method of
monotone iterations to deal with the existence of traveling wave fronts for
delayed equations (\ref{pde-0}) (see also \cite{ma,mazou,hualurua,huazou}).
This produces a monotone sequence of positive functions that converges to a
traveling wave front, which is an increasing positive solution of an equation
of the form
\begin{equation}
D\phi^{\prime\prime}(\xi)-c\phi^{\prime}(\xi)+f_{c}(\phi_{\xi})=0,\;\;\xi
\in\mathbb{R}.\label{w-0}%
\end{equation}
This method is contingent on the construction of a pair of upper and lower
solutions, which satisfy the inequality versions of (\ref{w-0}). It turns out
that this issue is closely related to the existence and uniqueness of smooth
and bounded solutions of nonhomogeneous equations of the form
\begin{equation}
Dx^{\prime\prime}(t)-cx^{\prime}(t)-\beta x(t)+g(t)=0,\quad\mbox
{for all}\ \ t\in\mathbb{R},\label{per}%
\end{equation}
where the constant $c,\beta>0$. In differential equations results of this type
are known as Perron Theorem, see \cite{chilat,hinnaiminshi,murnaimin,naimin}
and for interesting applications see \cite{chohal,chilat} and the references therein.

\medskip We would like to point out that the smoothness required for upper and
lower solutions of (\ref{w-0}) obstruct the search for bounded
solutions of (\ref{per}). One faces the following dilemma, see
e.g. \cite{zouwu,wuzou,ma}. An excessive relaxation of the
smoothness of upper and lower solutions of (\ref{w-0}) simplifies
their finding, but would not generate the sought monotone
iteration scheme. The main reason is due to the failure of Perron
Theorem for this class of weak solutions of (\ref{per}).

In the light of the above remark, the question of existence of traveling waves
front in predator-prey or Belousov- Zhabotinskii models with delay, as
addressed  in \cite{ma}, remains open.

The purpose of this paper is two folds: First to set up a rigorous framework
for the Monotone Iteration Method and then apply it to the predator-prey and
Belousov- Zhabotinskii models with delay.

\medskip We now briefly outline the plan of this paper. In the next section,
to recall the main concepts and tools, we discuss a modified version of Perron
Theorem for $C^{1}$-solutions of (\ref{per}). Next we show how crucial it is
to the monotone iteration method. Remarks and counter-examples are used to
explain the pitfalls of non-smooth upper solutions. In Theorem \ref{the 1} one
finds a rigorous   framework to construct fail-safe upper and lower solutions
which are then used for the delayed predator-prey and Belousov-Zhabotinskii
equations, see Theorems  \ref{the 2}, \ref{the 3}.

To conclude we would like to emphasise that Perron Thereom dictates $C^{1}%
$-smoothness which makes the search for upper and lower solutions much harder
than $C^{0}$-smooth solutions as in \cite{ma}, and of course easier than
$C^{2}$ solutions, which seems to be impossible for the above mentioned equations.

\section{The Monotone Iteration Scheme}

In this section we introduce a modified version of Perron Theorem which is
central to a rigorous framework for the monotone iteration scheme. We show
that this framework provides a clean procedure for constructing upper and
lower solutions of reaction diffusion equation with delay.

\subsection{Bounded solutions of nonhomogeneous equations}

The theory of bounded solutions of nonhomogeneous equations is a classical
topic of the theory of ordinary differential equations, which can be found in
\cite{chilat,hinnaiminshi,liunguminvu,murnaimin} and the references therein.

Let us first consider the equation
\begin{equation}
u^{\prime\prime}(t)+\alpha u^{\prime}(t)+\beta u(t)+f(t)=0,\quad
t\in\mathbb{R},\;u(t)\in\mathbb{R},\label{ode}%
\end{equation}
where $f$ is a function that is continuous and bounded on $\mathbb{R}%
\backslash\{0\}$ and has the right and left limits at $x=0$, $f(0^{+})$ and
$f(0^{-});$ We always assume that $\alpha\;$and $\beta$ are real numbers with
$\beta<0$ so that the characteristic equation
\[
\lambda^{2}+\alpha\lambda+\beta=0
\]
has two distinct real roots of opposite signs $\lambda_{1}<0$ and
$\lambda _{2}>0$. To our knowledge, all available results on
classical solutions deal with bounded and continuous forcing terms
$f$ on the entire real line (see e.g.
\cite{chilat,hinnaiminshi,liunguminvu,murnaimin}). As we allow
jump discontinuities in $f$ in (\ref{ode}), we need to modify the
concept of solutions as well as the conditions for their existence
and boundedness. We first need a definition.

\begin{definition}
Suppose that $f$ is a bounded and continuous function on $\mathbb{R}%
\backslash\{0\}$ and both $f(0^{+})$ and $f(0^{-})$ exist. Then, a function
$u$ defined on $\mathbb{R}$ is said to be a generalized solution of
(\ref{ode}) if

\begin{enumerate}
\item $u$ and $u^{\prime}$ are bounded and continuous on $\mathbb{R}$,

\item $u^{\prime\prime}$ exists and is continuous on $\mathbb{R}%
\backslash\{0\}$, and both $u^{\prime\prime}(0^{-})$ and $u^{\prime\prime
}(0^{+})$ exist.
\end{enumerate}
\end{definition}

Below is a version of Perron Theorem for generalized solutions with
discontinuous $f$.

\begin{lemma}\label{lem per}
\label{lem 1} Consider equation (\ref{ode}) with $\beta<0$, and
assume that

\begin{enumerate}
\item $f$ is a bounded and continuous function on
$\mathbb{R}\backslash\{0\}$ and both $f(0^{+})$ and $f(0^{-})$
exist,
\item  Equation (\ref{ode}) holds in the classical sense for
all $t$ except possibly at $t=0$.
\end{enumerate}

Then, Eq. (\ref{ode}) has a unique generalized solution $u$ given by
\begin{equation}
u(t)=Gf(t):=\frac{1}{\lambda_{1}-\lambda_{2}}\left(  \int_{-\infty}
^{t}e^{\lambda_{1}(t-s)}f(s)ds+\int_{t}^{+\infty}e^{\lambda_{2}(t-s)}%
f(s)ds\right)  .\label{bs}%
\end{equation}
\end{lemma}

\begin{proof}
Observe that function $Gf,$ which is defined by (\ref{bs}), exists
and is bounded. A simple computation shows that $Gf$ is in fact
continuously differentiable on $\mathbb{R}$, and
$(Gf)^{\prime\prime}$ exists and is continuous on the whole
$\mathbb{R}$ with a possible exception at $t=0$ at which both
$(Gf)^{\prime \prime}(0^{+})$ and $(Gf)^{\prime\prime}(0^{-})$
exist. It is important to observe that while $Gf$ and
$(Gf)^{\prime}$ are bounded on $\mathbb{R}$, $Gf$ is a particular
solution of (\ref{ode}) on each of the two disjoint intervals
$(-\infty,0)$ and $(0,+\infty)$. Obviously, $v=u-Gf$ is then a
classical solution of the homogeneous equation associated with
(\ref{ode}) on $(0,\infty),$ and hence
\[
u(t)-Gf(t)=ae^{\lambda_{1}t}+be^{\lambda_{2}t}\;\;\;\;\text{for}\;\;t>0
\]
where $a,b$ are constants. From the boundedness of $u$ and $Gf$ on
$(0,\infty)$ it follows that $b=0$ i.e.%
\[
u(t)-Gf(t)=ae^{\lambda_{1}t}\;\;\;\;\text{for}\;\;t>0.
\]
Similarly the classical solution $u-Gf$ on $(-\infty, 0)$ is of the
form
\[
u(t)-Gf(t)=de^{\lambda_{2}t}\;\;\;\;\;\text{for }t<0.
\]
Combining both behaviors we have
\begin{equation}
u(t)-Gf(t)= \left\{
\begin{array}
[c]{l}%
ae^{\lambda_{1}t},\quad t>0,\\
de^{\lambda_{2}t},\quad t\leq0.
\end{array}\right.
\label{ugf}
\end{equation}
Use the fact that although $u(t)-Gf$ has no second derivative at
$t=0$, it is still continuously
differentiable on the real line, and in particular at $t=0$, to deduce the interface conditions%
\[
\left\{
\begin{array}
[c]{c}%
a=d\\
\lambda_{1}a=\lambda_{2}d.
\end{array}
\right.
\]
Since $0\not=\lambda_{1}\not =\lambda_{2}\not= 0$, we must have
$a=d=0$, that is $u-Gf=0$. This completes the proof of the lemma.
\end{proof}

\begin{remark}
Lemma \ref{lem 1} shows that the continuity of $u^{\prime}$ plays
a crucial role in the uniqueness of the solution. The failure of
the continuity even at a single point, would lead to a
non-uniqueness of bounded solutions. We now illustrate this
important fact by a simple and yet explicit example.
\end{remark}

\begin{example}\label{ex cou}
\label{ex 1} Consider bounded solutions for the equation
\begin{equation}
y^{\prime\prime}(t)-y(t)=0\;\;\text{where
}t\in\mathbb{R.}\label{ex2}
\end{equation}
Since $f=0,$ by Lemma \ref{lem 1} the only bounded continuously
differentiable solution (that is, generalized solution) would be
$0.$ However, if we allow for a non-differentiable function to be a
solution, then the continuous function
\begin{equation}
y(t):=\left\{
\begin{array}
[c]{l}%
e^{-t},\quad t>0,\\
e^{t},\quad t\leq0 ,
\end{array}
\right.  \label{yex}%
\end{equation}
 satisfies (\ref{ex2}) for all
$t\in\mathbb{R}\backslash\{0\}$. Obviously, it is a nonzero bounded
and continuous function on $\mathbb{R}$, but not differentiable at
$t=0$.
\end{example}

\begin{remark}\label{rem 2}
The function $y$ in (\ref{yex}) can also serve as a counter-example
to the identity (see the proof of \cite[Lemma 3.5]{ma})
\begin{align}
&  \frac{d}{\lambda_{2}-\lambda_{1}}\left(
\int_{-\infty}^{t}e^{\lambda
_{1}(t-s)}\varphi(s)ds+\int_{t}^{\infty}e^{\lambda_{2}(t-s)}\varphi(s)\right)
\label{xuwu}\\
&
\;\;\;\;\;\;\;\;\;\;\;\;=y(t)+\frac{1}{\lambda_{2}-\lambda_{1}}\sum
_{j=k+1}^{m}e^{\lambda_{1}(t-T_{j})}(y^{\prime}(T_{j}+)-y^{\prime
}(T-))\nonumber\\
&  \;\;\;\;\;\;\;\;\;\;\;\;\;\;\;\;-\sum_{j=1}^{k}e^{\lambda_{2}(t-T_{j}%
)}(y^{\prime}(T_{j}+)-y^{\prime}(T-)),\nonumber
\end{align}
where $y$ is a piece-wise $C^{2}$ solution of
\[
dy^{\prime\prime}(t)-cy^{\prime}(t)-\beta y(t)=\varphi(t)\;\;\;\;\text{on\ \ }%
\mathbb{R}\backslash\{T_{0},...,T_{m}\}.
\]
In the particular case when $d=1,\;\lambda_{1}=-1,\;\lambda_{2}=1$,
$c=0,\;\beta=1$, $m=0$, $T_{0}=0,$ the expression (\ref{yex})
defines a solution to $y^{\prime\prime}(t)-y(t)=0$ for all $t\not
=0.$ Since $\varphi=0$, the left hand side of (\ref{xuwu}) is zero
while the right
hand side is not%
\begin{align*}
0 &
=y(t)+\frac{1}{2}\sum_{j=k+1}^{m}e^{-t}(y^{\prime}(0^{+})-y^{\prime
}(0^{-}))-\frac{1}{2}\sum_{j=1}^{k}e^{t}(y^{\prime}(0^{+})-y^{\prime}%
(0^{-}))\\
0 &  =y(t)+\frac{1}{2}\sum_{j=k+1}^{m}e^{-t}(-2)-\frac{1}{2}\sum_{j=1}%
^{k}e^{t}(-2)\\
0 &  =y(t)-e^{-t}+e^{t}\neq0
\end{align*}

We hope that that the argument for a continuous $y^{\prime}$ is by
now clear and self-evident.
\end{remark}

\subsection{Delayed reaction diffusion equations}

Consider the following system of reaction-diffusion equations with
time delay
\begin{equation}\label{pde}
\frac{\partial u(x,t)}{\partial t}=D\frac{\partial^{2}u(x,t)}{\partial x^{2}%
}+f(u_{t}),
\end{equation}
where $t\geq0$, $x\in\mathbb{R}$, $u(x,t)\in\mathbb{R}^{n}$, $D=diag(d_{1}%
,...,d_{n})$ with $d_{i}>0$ for all $i=1,...,n$, $f:C([-\tau,0],\mathbb{R}%
^{n})\rightarrow\mathbb{R}^{n}$ is a continuous, and $u_{t}(x)$ is
an element of $C([-\tau,0],\mathbb{R}^{n})$, defined as
\[
u_{t}(x)=u(t+\theta,x),\quad\theta\in\lbrack-\tau,0],\;t\geq0,\;x\in
\mathbb{R}.
\]
Throughout this paper we will assume that
\begin{equation}
f(\hat{0})=f(\hat{K})=0\ \mbox{and}\ f(\hat{u})\not =0,\ \mbox{for}%
\ u\in\mathbb{R}^{n}\ \mbox{with}\ \mathbf{0}<u<\mathbf{K},
\end{equation}
where $\hat{0}$ denotes the function $\varphi\in
C[-\tau,0],\mathbb{R}^{n})$ such that $\varphi(\theta)={\bf 0}:= \{
0,0,...,0\}^T \in {\mathbb R}^n$ for all $\theta\in\lbrack-\tau,0]$
and $\hat{K}$ denotes the function $\psi\in
C[-\tau,0],\mathbb{R}^{n})$ such that $\psi(\theta)={\bf
K}:=(K_1,K_2,...,K_n)^{T}\in\mathbb{R}^{n}$, with given positive
$K_i$, $i=1,2,...,n$. In order to use comparison arguments, we use
the natural partial order in
$\mathbb{R}^{n}$ to compare two vectors $x=(x_{1},...,x_{n})^{T}%
,\;y=(y_{1},...,y_{n})^{T}\in\mathbb{R}^{n}$, that is $x\geq y$ if
and only if $x_{i}\geq y_{i}$, for all $i=1,...,n$, and if there is
at least an $i$ in $\{ 1,...,n\}$ such that $x_i<y_i$, we write
$x<y$. The "interval" $[{\bf 0},{\bf K}]$ consists of all vectors
$v\in {\mathbb R}^n$ such that $0\le v_i \le K_i$ for all
$i=1,2,...,n$.

As usual we look for traveling wave solutions of (\ref{pde}) in the
form $u(x,t)=\phi(x+ct)$, where $\phi\in
C^{2}(\mathbb{R},\mathbb{R}^{n})$, and $c>0$ is a constant.
Substituting $u(x,t)=\phi(x+ct)$ into (\ref{pde}) leads to the
following \textit{wave equation}
\begin{equation}\label{wave}
D\phi^{\prime\prime}(\xi)-c\phi^{\prime}(\xi)+f_{c}(\phi_{\xi})=0,\;\;\xi
\in\mathbb{R},
\end{equation}
where
$f_{c}:X_{c}:=C([-c\tau,0],\mathbb{R}^{n})\rightarrow\mathbb{R}^{n}$,
defined as
\[
f_{c}(\psi)=f(\psi^{c}),\quad\psi^{c}(\theta):=\psi(c\theta),\quad\theta
\in\lbrack-\tau,0].
\]
Throughout this section we assume that there exists a matrix $\beta
=diag(\beta_{1},...,\beta_{n})$ with $\beta\geq0$ such that
\begin{equation}\label{A}
f_{c}(\phi)-f_{c}(\psi)+\beta\lbrack\phi(0)-\psi(0)]\geq\mathbf{0}%
,\;\;\mbox{for all}\;\;\ \phi,\psi\in X_{c},\quad\mathbf{0}\leq\psi
(\theta)\leq\phi(\theta)\leq\mathbf{K},
\end{equation}
where $\mathbf{0}:=(0,...,0)^{T}\in\mathbb{R}^{n}$ and $\mathbf{K}%
:=(K,...,K)^{T}\in\mathbb{R}^{n}$.

The main purpose of this paper is to look for solutions $\phi$ of
(\ref{wave}) in the following subset of
$C(\mathbb{R},\mathbb{R}^{n})$
\[
\Gamma:=\{\varphi\in
C(\mathbb{R},\mathbb{R}^{n}):\varphi\;\;\;\mbox{is
nondecreasing, and}\ \ \lim_{\xi\rightarrow-\infty}\varphi(\xi)=\mathbf{0}%
,\ \lim_{\xi\rightarrow+\infty}\varphi(\xi)=\mathbf{K}\}
\]
The solution $\phi$ is called a \textit{monotone wave front} of
(\ref{pde}). Next we define an operator $H:BC(\mathbb{R},\mathbb{R}%
^{n})\rightarrow BC(\mathbb{R},\mathbb{R}^{n})$ by
\begin{equation}\label{A2}
H(\phi)(t)=f_{c}(\phi_{t})+\beta\phi(t),\;\;\;\;\;\;\;\;\phi\in C(\mathbb{R}%
,\mathbb{R}^{n}).
\end{equation}
The following monotonicity lemma was proved in \cite{wuzou}.

\begin{lemma}
\label{lem 3.1} Assume that (\ref{A}) and (\ref{A2}) hold. Then, for
any $\phi\in\Gamma,$ we have that

\begin{enumerate}
\item $H(\phi)(t)\geq0,\ t\in\mathbb{R}$,

\item $H(\phi)(t)$ is nondecreasing in $t\in\mathbb{R,}$

\item $H(\psi)(t)\leq H(\phi)(t)$ for all $t\in\mathbb{R}$, if $\psi\in
C(\mathbb{R},\mathbb{R}^{n})$ is given so that
$\mathbf{0}\leq\psi(t)\leq \phi(t)\leq\mathbf{K}$ for all
$t\in\mathbb{R}$.
\end{enumerate}
\end{lemma}

Our definition of upper solutions, which is found below, requires
more smoothness and boundedness than those in \cite{wuzou} and
\cite{ma}.

\begin{definition}
A function $\rho\in C^{2}(\mathbb{R},\mathbb{R}^{n})$ with
$\rho,\rho^{\prime },\rho^{\prime\prime}$ being bounded on
$\mathbb{R}$ is said to be an \textit{upper solution} of
(\ref{wave}) if it satisfies the following
\begin{equation}
D\rho^{\prime\prime}(t)-c\rho^{\prime}(t)+f_{c}(\rho_{t})\leq\mathbf{0}%
,\ \mbox{for all}\ \ t\in\mathbb{R}.\label{usol}%
\end{equation}
\end{definition}

\begin{definition}
A function $\rho\in C^{1}(\mathbb{R},\mathbb{R}^{n})$ is said to be
a quasi-upper solution of (\ref{wave}) if

\begin{enumerate}
\item
\[
\sup_{t\in\mathbb{R}}\Vert\rho(t)\Vert<\infty,\quad
\sup_{t\in\mathbb{R}}\Vert \rho^{\prime}(t)\Vert<\infty,
\]

\item $\rho^{\prime\prime}(t)$ exists and is continuous on $\mathbb{R}%
\backslash\{0\}$, and
\[
\sup_{t\in\mathbb{R}\backslash\{0\}}\Vert\rho^{\prime\prime}(t)\Vert<
\infty,
\]

\item $\lim_{t\rightarrow0^{-}}\rho^{\prime\prime}(t)$ and $\lim
_{t\rightarrow0^{+}}\rho^{\prime\prime}(t)$ exist,

\item $\rho(t)$ satisfies
\begin{equation}
D\rho^{\prime\prime}(t)-c\rho^{\prime}(t)+f_{c}(\rho_{t})\leq\mathbf{0}%
,\ \mbox{for all}\ \ t\in\mathbb{R}\backslash\{0\}.\label{qsol}%
\end{equation}
\end{enumerate}
\end{definition}

Similarly, we define the concept of lower and quasi-lower solutions
to (\ref{wave}) by switching the inequalities in (\ref{usol}) and
(\ref{qsol}). Let $\phi\in C(\mathbb{R},\mathbb{R}^{n})$ such that
$\mathbf{0}\leq \phi(t)\leq\mathbf{K}$. By Lemma \ref{lem 3.1}, we
have $\mathbf{0}\leq H(\phi)(t)\leq\mathbf{K}$, so $H(\phi)$ is a
bounded function. Note that
\[
\lambda_{1i}:=\frac{c-\sqrt{c^{2}+4\beta_{i}d_{i}}}{2d_{i}}<0,\quad
\lambda_{2i}:=\frac{c+\sqrt{c^{2}+4\beta_{i}d_{i}}}{2d_{i}}>0
\]
are the real roots of the equation
\[
d_{i}\lambda^{2}-c\lambda-\beta_{i}=0,\;\;\;\;\;\;i=1,\;2,...,\;n.
\]
The existence of a unique generalized bounded solution of
\begin{equation}
Dx^{\prime\prime}(t)-cx^{\prime}(t)-\beta x(t)+H(\phi)(t)=0,\quad
\mbox{for all}\ \ t\in\mathbb{R}%
\end{equation}
follows by Perron Theorem, see Lemma \ref{lem 1}, is
\[
G(H(\phi))=(G_{1}(H_{1}(\phi),\;G_{2}(H_{2}(\phi)),...,\;G_{n}(H_{n}%
(\phi)))^{T}%
\]
where
\begin{equation}
G_{i}(H(\phi)):=\frac{1}{\lambda_{1i}-\lambda_{2i}}\left(  \int_{-\infty}%
^{t}e^{\lambda_{1i}(t-s)}H_{i}(\phi)(s)ds+\int_{t}^{+\infty}e^{\lambda
_{2i}(t-s)}H_{i}(\phi)(s)ds\right)  .
\end{equation}

The Monotone Iteration Scheme is constructed as follows: We start
out with a quasi-upper solution $\phi_{0}$, and then use the
recurrence formula
\begin{equation}\label{11}
\phi_{n}:=G(H(\phi_{n-1})),\quad n=1,2,....
\end{equation}

\begin{lemma}
\label{lem mono} Assume that $\phi,\phi_{0}\in\Gamma$ are,
respectively, a quasi-lower and quasi-upper solutions of
(\ref{wave}), such that $\phi (t)\leq\phi_{0}(t)$ for all
$t\in\mathbb{R}$. Then

\begin{enumerate}
\item[(i)]$\phi_{1}\in\Gamma$ for all $t\in\mathbb{R}$,

\item[(ii)] $\phi_{1}$ is a upper solution of (\ref{wave}), and
\[
\phi(t)\leq\phi_{1}(t)\leq\phi_{0}(t),\quad\mbox{for all}\
t\in\mathbb{R}.
\]
\end{enumerate}
\end{lemma}

\begin{proof}
For the proof of (i) we refer to \cite[Lemma 3.3]{wuzou}. \medskip
Next, we
prove (ii). Clearly $\phi_{1}=G(H(\phi_{0}))\in C^{2}(\mathbb{R}%
,\mathbb{R}^{n})$ and by Lemma \ref{lem 1}, $\phi_{1}$ is also the
unique  bounded (classical) solution of the equation
\[
Dx^{\prime\prime}(t)-cx^{\prime}(t)-\beta x(t)+H(\phi_{0})(t)={\bf
0} ,\quad\mbox{for all}\ \ t\in\mathbb{R}.
\]
By Lemma \ref{lem 3.1}, $H(\phi_{0})\geq0$, and so
\[
D\phi_{1}^{\prime\prime}(t)-c\phi_{1}^{\prime}(t)-\beta\phi_{1}(t)\leq
{\bf 0}\quad\mbox{for all}\ t\in\mathbb{R},
\]
which means that $\phi_{1}$ is a upper solution of (\ref{wave}) and
from its definition $\phi_{1}\geq\mathbf{0}$. We now show that the
sequence $\{\phi _{n}\}$ is decreasing. To this end, set
\begin{align}
w(t) &  :=\phi_{1}(t)-\phi_{0}(t),\nonumber\\
r(t) &  :=-Dw^{\prime\prime}(t)-cw^{\prime}(t)+\beta w(t),\quad t\in
\mathbb{R}.\label{r}%
\end{align}
By (\ref{11}) and the assumption that $\phi_{0}$ is a quasi-upper
solution of (\ref{wave}) we see that $r_{i}(t)$ is non-positive and
bounded for every $i=1,...,n$ and $t\in\mathbb{R}\backslash\{0\}$.
Moreover, $w_{i}(\cdot)$ and the $i$-component of (\ref{r})
satisfies all conditions of Lemma \ref{lem 1}. Therefore, it follows
that as a bounded generalized solution of (\ref{r}), $w_i$ satisfies
\[
w_{i}(t)=Gr_{i}(t):=\frac{1}{\lambda_{1}-\lambda_{2}}\left(
\int_{-\infty
}^{t}e^{\lambda_{1}(t-s)}r_{i}(s)ds+\int_{t}^{+\infty}e^{\lambda_{2}%
(t-s)}r_{i}(s)ds\right)  \leq0
\]
which yields $w(t)=\phi_{1}-\phi_{0}\leq {\bf 0}$,\ i.e.
$\phi_{1}\leq\phi_{0}$.
\medskip Next we show that $\phi\leq\phi_{1}$. To this end, we again set
\begin{align}
v(t) &  :=\phi_{1}(t)-\phi(t),\nonumber\\
s(t) &  :=-Dv^{\prime\prime}(t)-cv^{\prime}(t)+\beta v(t)\geq
0,\;\;\;\;\;\;t\in\mathbb{R}.\label{s}%
\end{align}
By the definition of quasi-lower solution and $\phi_{1}$ and by
Lemma \ref{lem 1} we have
\begin{equation}
v_{i}(t)=G\left(  s_{i}\right)
(t):=\frac{1}{\lambda_{1}-\lambda_{2}}\left(
\int_{-\infty}^{t}e^{\lambda_{1}(t-\xi)}s_{i}(\xi)d\xi+\int_{t}^{+\infty
}e^{\lambda_{2}(t-\xi)}v_{i}(\xi)d\xi\right)  \geq0\label{w}%
\end{equation}
which implies $\phi_{1}(t)\geq\phi(t)$. $\Box$
\end{proof}

\begin{remark}
As shown in Example \ref{ex 1}, without assumption on the continuity
of the derivative of $\phi_{0}$, formula (\ref{w}) may not be true
because we know only that $w$ is bounded and it may not be of $C^2$.
\end{remark}

Below we will use the notation $BC_{[\mathbf{0},\mathbf{K}
]}(\mathbb{R},\mathbb{R}^{n}):= \{ g\in
BC(\mathbb{R},\mathbb{R}^{n})|\ {\bf 0} \le g(t) \le {\bf K}, \
\forall t\in {\mathbb R}\} $ as a closed convex subspace of the
Banach space $BC(\mathbb{R},\mathbb{R}^{n})$ equipped with the
sup-norm. We are now ready to state the Monotone Iteration Method.

\begin{theorem}\label{the 1} Assume that

\begin{enumerate}
\item  there exist a quasi-upper solution of (\ref{wave}) $\overline{\phi}%
_{0}\in\Gamma$ and a non-zero non-decreasing quasi-lower solution
$\underline{\phi_{0}}$ such that
$0<\underline{\phi_{0}}(t)\leq\overline {\phi_{0}}(t)$ for all
$t\in\mathbb{R}$,

\item $H$ is continuous on $BC_{[\mathbf{0},\mathbf{K}]}(\mathbb{R}%
,\mathbb{R}^{n})$, and
\begin{equation}
\sup_{\varphi\in\Gamma}\Vert H(\varphi)\Vert<\infty.
\end{equation}
\end{enumerate}

Then, the following assertions hold:

\begin{enumerate}
\item [(i)]The sequence $\{\phi_{n}\}_{n=1}^{\infty}$, defined as
above, is a decreasing sequence in $\Gamma$,

\item[(ii)] $\lim_{n\rightarrow\infty}\phi_{n}(t)=\phi(t)$ is a
monotone wave front of (\ref{wave}). Moreover, this limit is uniform
on each compact interval of the real line.
\end{enumerate}
\end{theorem}

\begin{proof}
By Lemma \ref{lem mono} Claim (i) is clear. To prove (ii) we use
Arzela-Ascoli
Theorem. Note that by Lemma \ref{lem mono}, all $\{\phi_{n}\}_{n=1}^{\infty}%
$\ $\;$are $C^{2}$ functions. From
\begin{align*}
(\phi_{n+1})_{i}(t) &  =G_{i}(H(\phi_{n}))(t)\\
&  =\frac{1}{\lambda_{1i}-\lambda_{2i}}\left(  \int_{-\infty}^{t}%
e^{\lambda_{1i}(t-s)}H_{i}(\phi_{n})(s)ds+\int_{t}^{+\infty}e^{\lambda
_{2i}(t-s)}H_{i}(\phi_{n})(s)ds\right)  ,
\end{align*}
we deduce
\begin{align}
|(\phi_{n+1})_{i}^{\prime}(t)| &
\leq\frac{2\sup_{\varphi\in\Gamma}\Vert H(\varphi)\Vert}{|\lambda
_{1i}-\lambda_{2i}|}.
\end{align}
Therefore, on each interval $[-N,\;N],$ the set
$\{\phi_{n}\}_{n=1}^{\infty}$ is equicontinuous and by Arzela-Ascoli
Theorem, we have a uniformly convergent subsequence
$\{\phi_{n_{k}}\}_{k}$ on $[-N,\;N].$ Since the sequence is monotone
we can conclude that the sequence $\{\phi_{n}\}_{n=1}^{\infty}$
converges uniformly to a continuous and non-decreasing function
$\phi$ such that $\mathbf{0}\leq\phi(t)\leq\mathbf{K}$ on every
compact interval of the
real line. Since $\phi\in BC_{[\mathbf{0},\mathbf{K}]}(\mathbb{R}%
,\mathbb{R}^{n})$ and $H$ is continuous on $BC_{[\mathbf{0},\mathbf{K}%
]}(\mathbb{R},\mathbb{R}^{n})$, for every fixed $(t,s)$, we have
\[
\lim_{n\rightarrow\infty}e^{\lambda_{1i}(t-s)}H_{i}(\phi_{n})(s)=e^{\lambda
_{1i}(t-s)}H_{i}(\phi)(s).
\]
Lebesgue's Dominated Convergence Theorem, then yields
\begin{align*}
(\phi(t))_{i} &  =\lim_{n\rightarrow\infty}(\phi_{n+1}(t))_{i}
 =\frac{1}{\lambda_{1i}-\lambda_{2i}}\left(  \int_{-\infty}^{t}%
e^{\lambda_{1i}(t-s)}H_{i}(\phi)(s)ds+\int_{t}^{+\infty}e^{\lambda_{2i}%
(t-s)}H_{i}(\phi)(s)ds\right)  \\
&  =(G(H(\phi))(t))_{i}.
\end{align*}
This shows in particular that $\phi$ is a classical solution of the equation%
\[
D\phi^{\prime\prime}(t)-c\phi^{\prime}(t)-\beta\phi(t)+H(\phi)(t)={\bf
0}.
\]
It remains to see that $\phi\in\Gamma$. First since
$\phi(t)\leq\phi_{n}(t)$
for all $n$ and $t\in\mathbb{R}$%
\[
\mathbf{0}\leq\limsup_{t\rightarrow-\infty}\phi(t)\leq\limsup_{t\rightarrow
-\infty}\phi_{n}(t)=\mathbf{0}.
\]
Next to show that $\lim_{t\rightarrow\infty}\phi(t)=\mathbf{K,}$ we
use $\mathbf{0}\not =\underline{\phi_{0}}(t)\leq\phi(t)$ for all
$t$, and the fact that $\phi$ is non-decreasing means that
$\lim_{t\rightarrow\infty}\phi(t)$ exists as a vector $\mathbf{L}\in
[\mathbf{0},\mathbf{K}]$. Obviously, $\mathbf{L}\not =\mathbf{0}$,
and so $\mathbf{0}<\mathbf{L}$. We now use the
delay and as in \cite[Proposition 2.1]{wuzou}, we must have $f_{c}%
(\mathbf{L})=0$. But the assumption on the function $f$ makes it
impossible if $\mathbf{L}\not =\mathbf{K}$. So,
$\lim_{t\rightarrow\infty}\phi (t)=\mathbf{K}$, and $\phi\in\Gamma$.
\end{proof}

\section{Traveling Waves for the Predator-Prey Model}

Consider the predator-prey model with diffusion and delay in time
\begin{equation}
\left\{
\begin{array}
[c]{l}%
\frac{\partial u(x,t)}{\partial
t}=d_{1}\frac{\partial^{2}u(x,t)}{\partial
x^{2}}+ru(x,t)\left[  \left(  1-\frac{1}{P}u(x,t)\right)  -av(x,t)\right] \\
\frac{\partial v(x,t)}{\partial
t}=d_{2}\frac{\partial^{2}v(x,t)}{\partial x^{2}}+v(x,t)\left[
-\nu+bu(x,t-\tau)\right]  ,
\end{array}
\right.  \label{pp}%
\end{equation}
where $u(x,t),v(x,t)$ are scalar functions,
$(x,t)\in\mathbb{R}\times [ \tau   ,\infty)$, and $d_{1}
,d_{2},r,P,a,b,\nu,\tau$ are all positive constants. Define
$\mathbb{R}$-valued functions $f_{1}(\phi,\psi)$ and
$f_{2}(\phi,\psi)$ as follows: For each $\phi,\psi\in
C([-\tau,0],\mathbb{R})$, the functionals
\begin{align*}
f_{1}(\phi,\psi)  &  =r\phi(0)\left[  \left(
1-\frac{\phi(0)}{P}-a\psi
(0)\right)  \right] \\
f_{2}(\phi,\psi)  &  =\psi(0)\left(  -\nu+b\phi(-\tau )\right)  .
\end{align*}
If we assume that
\[
1>\frac{\nu}{Pb}\quad\Leftrightarrow\quad P>\frac{\nu}{b},
\]
then the model has a non trivial positive steady state
\[
\left(  \frac{\nu}{b},\;\frac{1}{a}(1-\frac{\nu}{Pb})\right)^{T}.
\]

The above model with $\tau =0$ has been treated in
\cite{mur,dun1,dun2} (for related results see
\cite{mis,mishut,misrei}.) The case with delay $\tau
>0$ was first considered in \cite{ma}. As noted in Remark \ref{rem
2}, the smoothness was overlooked in the proof of \cite[Lemma
3.5]{ma}. This gap is due to the failure of Perron Theorem for
this type of "suppersolutions". Hence, the question on the
existence of traveling waves for this model with delay remains
open. A quick remedy is provided by the results of the previous
section, which require us to construct  quasi-upper and
quasi-lower solutions which have more smoothness than
"suppersolutions" in \cite{ma}.

\subsection{Wave equations and upper, lower solutions}

Let us agree to denote by $\varphi_{1}(x+ct):=u(x,t)$, $\varphi_{2}%
(x+ct):=v(x,t),$ and $\xi=x+ct$ .  We recast (\ref{pp}) into
\[%
\begin{array}
[c]{l}%
c\varphi_{1}^{\prime}(\xi)=d_{1}\varphi_{1}^{\prime\prime}(\xi)+r\varphi
_{1}(\xi)\left[  \left(  1-\frac{\varphi_{1}(\xi)}{P}\right)
-a\varphi
_{2}(\xi)\right]  =0\\
c\varphi_{2}^{\prime}(\xi)=d_{2}\varphi_{2}^{\prime\prime}(\xi)+\varphi
_{2}(\xi)\left[  -\nu+b\varphi_{1}(\xi-c\tau)\right]  =0,
\end{array}
\]
which can be re-written as
\begin{equation}%
\begin{array}
[c]{l}%
d_{1}\varphi_{1}^{\prime\prime}(\xi)-c\varphi_{1}^{\prime}(\xi)+r\varphi
_{1}(\xi)\left[  \left(  1-\frac{\varphi_{1}(\xi)}{P}\right)
-a\varphi
_{2}(\xi)\right]  =0\\
d_{2}\varphi_{2}^{\prime\prime}(\xi)-c\varphi_{2}^{\prime}(\xi)+\varphi
_{2}(\xi)\left[  -\nu+b\varphi_{1}(\xi-c\tau)\right]  =0.
\end{array}
\label{we}%
\end{equation}
We will find monotone solutions of the above wave equation such that
\begin{align}\label{boundary}
\lim_{t\rightarrow-\infty}\varphi_{1}(t)  &
=0\;\;\;\;\;\;\lim_{t\rightarrow
+\infty}\varphi_{1}(t)=p:=\frac{\nu}{b}\\
\lim_{t\rightarrow-\infty}\varphi_{2}(t)  &  =0\;\;\;\;\;\;\;\lim
_{t\rightarrow+\infty}\varphi_{2}(t)=q:=\frac{1}{a}\left(
1-\frac{\nu} {Pb}\right)  .
\end{align}

We now check the quasi-monotonicity of $(f_1,f_2)^T$, that is, there
are positive numbers $\beta_{1},\beta_{2}$ such that
\begin{align*}
f_{1}(\phi_1,\psi_2)-f_1(\psi_1,\psi_2)+\beta_{1}(\phi (0) -\psi (0))  &  \ge 0\\
f_2(\phi_1,\phi_2) -f_{2}(\psi_1,\psi_2) + \beta (\phi_2(0)-\psi_2
(0)) & \ge 0
\end{align*}
for all
$${\bf 0} \le \psi (s) \le \phi (s) \le {\bf K} :=
\left(\frac{\nu}{b}, \frac{1}{a}\left(1-\frac{\nu}{Pb}\right)
\right)^T.$$ Using the definition of the functions $(f_1,f_2)^T$ we
can easily come up with a pair of positive numbers $\beta_1,
\beta_2$ that satisfies
$$
\begin{cases}
\beta_1 \ge \max \left( -r +\frac{r\nu}{Pb} +\frac{\nu}{b}, 0
\right) \\
\beta_2 \ge \nu .
\end{cases}
$$

\subsection{Upper solutions}

Let us define $\varphi(t):=(\varphi_{1}(t),\varphi_{2}(t))^{T}$,
where
\begin{equation} \label{up}
\varphi_{1}(t)=\left\{
\begin{array}
[c]{l}%
\frac{\nu}{2b}e^{\lambda_{1}t},\ \ \;\;\;\;\;t\leq0\\
\frac{\nu}{b}-\frac{\nu}{2b}e^{-\lambda_{1}t},\ \ t>0
\end{array}
\right.  \;\;\;\;\text{and}\;\;\varphi_{2}(t)=\left\{
\begin{array}
[c]{l}%
\frac{q}{2}e^{\lambda_{2}t},\ \;\;\;\;\;\ t\leq0\\
q-\frac{q}{2}e^{-\lambda_{2}t},\ \ t>0,
\end{array}
\right.
\end{equation}
where the parameters are chosen so that $$
 \lambda_{1}  =\frac{c+\sqrt{c^{2}-4d_{1}r}}{2d_{1}}, \ \   \lambda_{2}  =\frac{c}{2d_{2}},     \ \
   p     =\frac{\nu}{b}, \ \   q =\frac{1}{a}\left( 1-\frac{\nu}{Pb}\right) .
$$
Therefore, $\lambda_{1}$ and $\lambda_{2}$ are simply the roots of
the quadratics
\begin{equation}
d_{1}\lambda^{2}-c\lambda+r=0,\;\;\;\;d_{2}\lambda^{2}-c\lambda=0. \label{l2}%
\end{equation}

\begin{lemma}\label{lem 101}
\label{lem up} Assume that $\nu<b$, and $d_{2}>2d_{1}$. Then, there
exists a constant $c^{\ast}=c^{\ast}(a,b,r,\nu,P)>0$ such that if
$c>c^{\ast}$,
$\overline{\phi_{0}}(t):=(\varphi_{1}(t),\varphi_{2}(t))^{T}$ is a
monotone quasi-upper solution of the wave equation.
\end{lemma}
\begin{proof} It is easily seen that both $\varphi_{1}$ and
$\varphi_{2}$ are continuously differentiable, monotone and
obviously
\begin{align*}
\lim_{t\rightarrow-\infty}\varphi_{1}(t)  &
=0\;\;\;\;\;\;\lim_{t\rightarrow
+\infty}\varphi_{1}(t)=p:=\frac{\nu}{b}\\
\lim_{t\rightarrow-\infty}\varphi_{2}(t)  &  =0\;\;\;\;\;\;\;\lim
_{t\rightarrow+\infty}\varphi_{2}(t)=q:=\frac{1}{a}\left(  1-\frac{\nu}%
{Pb}\right)  .
\end{align*}
Next, we have%

\begin{align*}
\varphi_{1}^{\prime}(t)  &  =\left\{
\begin{array}
[c]{l}%
\frac{p\lambda_{1}}{2}e^{\lambda_{1}t},\ \ t\leq0\\
\frac{p\lambda_{1}}{2}e^{-\lambda_{1}t},\ \ t>0
\end{array}
\right.  \quad\varphi_{1}^{\prime\prime}(t)=\left\{
\begin{array}
[c]{l}%
\frac{p\lambda_{1}^{2}}{2}e^{\lambda_{1}t},\ \ t\leq0\\
\frac{-p\lambda_{1}^{2}}{2}e^{-\lambda_{1}t},\ \ t>0
\end{array}
\right. \\
\varphi_{2}^{\prime}(t)  &  =\left\{
\begin{array}
[c]{l}%
\frac{q\lambda_{2}}{2}e^{\lambda_{2}t},\ \ t\leq0\\
\frac{q\lambda_{2}}{2}e^{-\lambda_{2}t},\ \ t>0
\end{array}
\right.  \quad\varphi_{2}^{\prime\prime}(t)=\left\{
\begin{array}
[c]{l}%
\frac{q\lambda_{2}^{2}}{2}e^{\lambda_{2}t},\ \ t\leq0\\
\frac{-q\lambda_{2}^{2}}{2}e^{-\lambda_{2}t},\ \ t>0
\end{array}
\right.  .
\end{align*}
A crucial property of the functions $\varphi_{1},\;\varphi_{2}$ is
their smoothness. We now look for conditions on the parameters so
that they form upper solutions. In doing so, we examine the cases
when $t\leq0$ and $t>0$ separately. Substituting the above
expressions into the first equation of (\ref{we}) and using
(\ref{l2}) yields for all $t\leq 0$
\begin{align*}
&  d_{1}\frac{p\lambda_{1}^{2}}{2}e^{\lambda_{1}t}-c\frac{p\lambda_{1}}%
{2}e^{\lambda_{1}t}+r\frac{p}{2}e^{\lambda_{1}t}-\frac{rp}{2P}e^{\lambda_{1}%
t}-\frac{arpq}{2}e^{\lambda_{1}t}e^{\lambda_{2}t}\\
&  =-\frac{rp}{2P}e^{\lambda_{1}t}-\frac{arpq}{2}e^{\lambda_{1}t}%
e^{\lambda_{2}t}\\
&  \leq0.\quad
\end{align*}
Similarly for $t>0$, we have
\begin{align*}
&  -d_{1}\frac{p\lambda_{1}^{2}}{2}e^{-\lambda_{1}t}-c\frac{p\lambda_{1}}%
{2}e^{-\lambda_{1}t}+r\left(  p-\frac{p}{2}e^{-\lambda_{1}t}\right)
-\frac
{r}{P}\left(  p-\frac{p}{2}e^{-\lambda_{1}t}\right) \\
&  \hspace{2cm}-ar\left(  p-\frac{p}{2}e^{-\lambda_{1}t}\right)
\left(
q-\frac{q}{2}e^{-\lambda_{2}(t-c\tau)}\right) \\
&  =\left(  rp-\frac{rp}{P}-{arpq}\right)  +\left(
\frac{rp}{2P}+\frac
{arpq}{2}-cp\lambda_{1}\right)  e^{-\lambda_{1}t}+\frac{arpq}{2}%
e^{-\lambda_{2}t}-\frac{arpq}{4}e^{-(\lambda_{1}+\lambda_{2})t} .
\end{align*}

To check the sign of the above expression, we first factor out $rp$,
which is positive,
\begin{equation}
\left(  1-\frac{1}{P}-{aq}\right)  +\left(  \frac{1}{2P}+\frac{aq}{2}%
-c\frac{\lambda_{1}}{r}\right)
e^{-\lambda_{1}t}+\frac{aq}{2}e^{-\lambda
_{2}t}-\frac{aq}{4}e^{-(\lambda_{1}+\lambda_{2})t}\label{up2} .
\end{equation}

A  sufficient condition on $P,$ can be obtained from (\ref{up2}) and
use if $\lambda_{2}>\lambda_{1}$,  which holds if  $d_{2}>2d_{1}$
\begin{align*}
& 2\left(  1-\frac{1}{P}-{aq}\right)  +\left(  \frac{1}{P}+aq-c\frac
{2\lambda_{1}}{r}\right)  e^{-\lambda_{1}t}+aqe^{-\lambda_{2}t}-\frac{1}%
{2}aqe^{-(\lambda_{1}+\lambda_{2})t}\\
& =\frac{2}{P}\left(  \frac{\nu}{b}-1\right)  +\left(  \frac{1}{P}%
+2aq-2c\frac{\lambda_{1}}{r}\right)  e^{-\lambda_{1}t} \leq 0
\end{align*}
to be non-positive for all $t>0,$ is%
\begin{align}
\frac{\nu}{b} &  <1\label{1}\\
\frac{1}{P}+2aq &  <2c\frac{\lambda_{1}}{r}\label{2}%
\end{align}
For (\ref{2}) to hold it is sufficient to take $c$ large enough.
More precisely, if  $c>\sqrt{4d_{1}r}$, then $\lambda_{1}>2c$ which,
in turn, leads to $c>\frac{1}{2}\sqrt{\left( \frac{1}{P}+2aq\right)
r.}$ Thus  the value of
$c^{\ast}$ can be estimated by%
\[
c^{\ast}=\max\left(  \frac{1}{2}\sqrt{\left(  \frac{1}{P}+2aq\right)
r},\;\sqrt{4d_{1}r}\right)  .
\]

\begin{corollary}
In case  $\frac{\nu}{b}<P<\frac{1}{2},$ then there $c^{\ast}=0.$
\end{corollary}
\begin{proof} It is enough to see that (\ref{up2}) can also be written as
\begin{align}
&  =\left(  1-\frac{1}{P}\right)  +\left(  \frac{1}{2P}-c\frac{\lambda_{1}}%
{r}\right)  e^{-\lambda_{1}t}+aq\left(
\frac{e^{-\lambda_{2}t}+e^{-\lambda
_{1}t}}{2}-1-\frac{e^{-(\lambda_{1}+\lambda_{2})t}}{4}\right)  \nonumber\\
&  =1-\frac{1}{2P}-c\frac{\lambda_{1}}{r}e^{-\lambda_{1}t}+\left(
\frac{e^{-\lambda_{1}t}-1}{2P}\right)  +aq\left(  \frac{e^{-\lambda_{2}%
t}+e^{-\lambda_{1}t}}{2}-1-\frac{e^{-(\lambda_{1}+\lambda_{2})t}}{4}\right)
  \leq 0.
\end{align}
It turns out that we have an upper solution, provided $c>0$ and
\begin{equation}
1-\frac{1}{2P}<0,\;\;\text{i.e. }\frac{\nu}{b}<P<\frac{1}{2}\text{\ }%
\label{coup1}
\end{equation}
no restrictions on $d_{1},$ $d_{2}$. \end{proof}

Now we  check the second equation of the wave equation. Starting
again with $t\leq0$ and substituting the expressions for
$\varphi_{1},\;\varphi_{2}$ into the second equation, together with
(\ref{l2}) yields
\begin{align*}
   d_{2}\frac{q\lambda_{2}^{2}}{2}e^{\lambda_{2}t}-c\frac{q\lambda_{2}}%
{2}e^{\lambda_{2}t}-\frac{q}{2}e^{\lambda_{2}t}\left(  \nu-b\frac{\nu}%
{2b}e^{\lambda_{1}(t-c\tau)}\right)
=-\frac{q\nu}{4}e^{\lambda_{2}t}\left(
2-e^{\lambda_{1}(t-c\tau)}\right)
 \leq0,\quad
\end{align*}
for all $ t\leq 0$ since $t-c\tau<0.$

Next, for $t>0$, substituting into the second equation of the wave
equation yields
\begin{align*}
&  -d_{2}\frac{q\lambda_{2}^{2}}{2}e^{-\lambda_{2}t}-c\frac{q\lambda_{2}}%
{2}e^{-\lambda_{2}t}-\nu\left(
q-\frac{q}{2}e^{-\lambda_{2}t}\right)
+b\left(  q-\frac{q}{2}e^{-\lambda_{2}t}\right)  \left(  p-\frac{p}%
{2}e^{-\lambda_{1}(t-c\tau)}\right)  \\
&  =-\left(  \nu-pb\right)  \frac{q}{2}+\left(  \nu-pb\right)  \frac{q}%
{2}\left(  e^{-\lambda_{2}t}-1\right)  -2c{\lambda_{2}}qe^{-\lambda_{2}%
t}-\frac{pqb}{2}e^{-\lambda_{1}(t-c\tau)}+\frac{pqb}{4}e^{-\lambda_{1}%
(t-c\tau)}e^{-\lambda_{2}t}.
\end{align*}
Since
\[
pb=\nu
\]
the above expression reduces to
\begin{align*}
&
-2c{\lambda_{2}}qe^{-\lambda_{2}t}-\frac{pqb}{2}e^{-\lambda_{1}(t-c\tau
)}+\frac{pqb}{4}e^{-\lambda_{1}(t-c\tau)}e^{-\lambda_{2}t}\\
&  \leq-2c{\lambda_{2}}qe^{-\lambda_{2}t}-\frac{q\nu}{2}e^{-\lambda
_{1}(t-c\tau)}+\frac{q\nu}{4}e^{-\lambda_{1}(t-c\tau)}e^{-\lambda_{2}t}\\
&
=-2cq\lambda_{2}e^{-\lambda_{2}t}-\frac{q\nu}{2}e^{-\lambda_{1}(t-c\tau
)}\left(  1\ -\frac{1}{2}e^{-\lambda_{2}t}\right)  \ \\
&  <0.
\end{align*}
and therefore the lemma is proved.\end{proof}

\subsection{Lower solutions}

We will construct a quasi-lower solutions to (\ref{we}) as follows:
Set
\[
\nu_{1}:=\frac{c}{d_{1}}.
\]
which is a root of
\[
d_{1}\nu_{1}^{2}-c\nu_{1}=0.
\]
Let us define functions $\psi_{2}(t)=0$ for all $t\in\mathbb{R}$,
and
\begin{equation}\label{psi1}
\psi_{1}(t)=\left\{
\begin{array}
[c]{l}%
\alpha ke^{\nu_{1}t},\;\;\;\;\;\;\;\;\;\;t<0\\
k-\alpha ke^{-\nu_{1}t},\;\;\;\;\;t\geq0,
\end{array}
\right.
\end{equation}
where the constants $k,$ and $\alpha$  are positive, sufficiently
small and will be determined at the end. A simple computation yields
\[
\psi_{1}^{\prime}(t):=\left\{
\begin{array}
[c]{l}%
\nu_{1}\alpha ke^{\nu_{1}t},\quad t<0\\
\nu_{1}\alpha ke^{-\nu_{1}t},\quad t\geq0.
\end{array}
\right.  ,\quad\psi_{1}^{\prime\prime}(t):=\left\{
\begin{array}
[c]{l}%
\nu_{1}^{2}\alpha ke^{\nu_{1}t},\quad t<0\\
-\nu_{1}^{2}\alpha ke^{-\nu_{1}t},\quad t\geq0.
\end{array}
\right.
\]
and the first equation of (\ref{we}) with $t\leq0$ gives
\begin{align*}
&  (d_{1}\alpha k\nu_{1}^{2}-c\alpha k\nu_{1})e^{\nu_{1}t}+r\alpha
ke^{\nu
_{1}t}\left[  1-\frac{1}{P}\alpha ke^{\nu_{1}t}\right]  \\
&  =r\alpha ke^{\nu_{1}t}\left[  1-\frac{1}{P}\alpha ke^{\nu_{1}t}\right]  \\
&  \geq0,
\end{align*}
provided $0<\alpha k<P$.

For $t>0$, substituting the derivatives of $\psi_{1}$ into the first
equation of (\ref{we}) leads to
\begin{align*}
&  -d_{1}\nu_{1}^{2}\alpha ke^{-\nu_{1}t}-c\nu_{1}\alpha ke^{-\nu_{1}%
t}+r(k-\alpha ke^{-\nu_{1}t})\left[  1-\frac{1}{P}(k-\alpha ke^{-\nu_{1}%
t})\right]  \\
&  =rk\left(  1-\frac{k}{P}\right)  +\alpha k\left(
\frac{2rk}{P}-2c\nu
_{1}-r\right)  e^{-\nu_{1}t}-\frac{r\alpha^{2}k^{2}}{P}e^{-2\nu_{1}t}\\
&  \geq 0.
\end{align*}
provided that $c$ and $k$ are fixed, $0<k<P$, and $\alpha$ is
sufficiently small.

Next, we need to choose $\alpha,\;k$ so that
\[
\psi_{1}(t)\leq\phi_{1}(t).
\]
In fact, a simple computation shows that we may choose $k\le
\nu/(4b)$ and $\alpha<(4Pb)/\nu$.

Obviously, this function $(\psi_{1},\psi_{2})^{T}$ is a quasi-lower
solution to (\ref{we}). Using the same procedure as for upper
solutions we will obtain a smooth lower solution as desired.

Therefore, we have proved

\begin{lemma}
Let all assumptions of Lemma \ref{lem 101} hold. Then, for every
fixed $c>c^*$, where $c^*$ is as in Lemma \ref{lem 101}, there are
(sufficiently small) constants $\alpha$ and $k$ such that the
function $\underline{\phi_{0}}(t):=(\psi_{1}(t),0)^{T}$, where
$\psi_{1}(t)$ is defined by (\ref{psi1}) is a quasi-lower solution
of (\ref{we}) that satisfies:
\begin{equation}
0<\underline{\phi_{0}}(t)\leq\overline{\phi_{0}}(t),\quad
t\in\mathbb{R},
\end{equation}
where $\overline{\ \phi_{0}}$ is the quasi-upper solution that is
defined by (\ref{up}) and mentioned in Lemma \ref{lem up}.
\end{lemma}

Finally, by Theorem \ref{the 1} we have
\begin{theorem}\label{the 2}
Assume that $\nu<b$, and $d_{2}>2d_{1}$. Then, there exists a
constant $c^{\ast}=c^{\ast}(a,b,r,\nu,P)>0$ such that if
$c>c^{\ast}$, Eq. (\ref{pp}) has a wave front solution
$u(x,t)=\varphi_1(x+ct)$, $v(x,t)=\varphi_2(x+ct)$. Moreover, if
$\frac{\nu}{b}<P<\frac{1}{2},$ $c*$ may be chosen to be $0$.
\end{theorem}

\section{Belousov- Zhabotinskii Equations}
In this section we will apply the results in Section 2 to prove the
existence of traveling waves to Belousov- Zhabotinskii Equations
\begin{equation}\label{pde-bz}
\begin{cases}
\frac{\partial }{\partial t}u(x,t) = \frac{\partial^2}{\partial
x^2} u(x,t) +u(x,t) [ 1-u(x,t) -r v(x,t-\tau )] ;\\
\frac{\partial }{\partial t}v(x,t) = \frac{\partial^2}{\partial x^2}
v(x,t) -bu(x,t)v(x,t) ,
\end{cases}
\end{equation}
where $r>0$, $b>0$ are constants, $u$ and $v$ are scalar functions.
We refer the reader to \cite{mur} for more details on the history as
well as applications of this kind of equations. Traveling waves in
these models without delay were considered in various papers (see
e.g. \cite{kan,kap,tro,yewan}). In the recent papers
\cite{ma,wuzou}, the models with delay were first studied as an
application of the Monotone Iteration Method, which was extended by
Wu and Zou. As we have noted earlier, in the previous model, the
Monotone Iteration Method reduces to finding a pair of upper and
lower solutions to the corresponding wave equation.  Again due to
the failure of Perron Theorem for the concepts of "upper solutions"
and "supersolutions" in \cite{wuzou,ma}, the question of existence
of traveling waves in these models remained open. Below, Perron
Theorem and Theorem \ref{the 1} will be used to answer this question
by constructing suitable quasi-upper and quasi-lower solutions to
the wave equation.

\medskip
The wave equation associated with Belousov- Zhabotinskii Equations
can be modified (see also \cite{wuzou,ma}) so that it takes the
form
\begin{equation}
\left\{
\begin{array}
[c]{l}%
\varphi_{1}^{\prime\prime}(t)-c\varphi_{1}^{\prime}(t)+\varphi_{1}(t)\left(
(1-r)-\varphi_{1}(t)+r\varphi_{2}(t-c\tau)\right)  =0\\
\varphi_{2}^{\prime\prime}(t)-c\varphi_{2}^{\prime}(t)+b\varphi_{1}(t)\left(
1-\varphi_{2}(t)\right)  =0.
\end{array}
\right.  \label{bz-we-0}%
\end{equation}
We  seek monotone solutions $(\varphi_{1},\varphi_{2})^{T}$ such
that
\begin{align*}
\lim_{t\rightarrow-\infty}\varphi_{1}(t)=0,\quad &
\lim_{t\rightarrow+\infty
}\varphi_{1}(t)=1\\
\lim_{t\rightarrow-\infty}\varphi_{2}(t)=0,\quad &
\lim_{t\rightarrow+\infty }\varphi_{2}(t)=1.
\end{align*}
To this end we recast the wave equations in the following form
\begin{align}
\lbrack\varphi_{1}^{\prime\prime}(t)-c\varphi_{1}^{\prime}(t)+\varphi
_{1}(t)]-r\varphi_{1}(t)(1-\varphi_{2}(t-c\tau))-\varphi_{1}^{2}(t)
&
=0\label{bz-we}\\
\left[
\varphi_{2}^{\prime\prime}(t)-c\varphi_{2}^{\prime}(t)\right]
+b\varphi_{1}(t)-b\varphi_{1}(t)\varphi_{2}(t) &  =0.\nonumber
\end{align}
Define the numbers $\lambda_{1}$ and $\mu_{1}$ as
\[
\lambda_{1}=\frac{c+\sqrt{c^{2}-4}}{2},\;\;\;\mu_{1}=\frac{c+\sqrt{c^{2}-4b}%
}{2},
\]
which are the  roots of the characteristic equations
\[
\lambda^{2}-c\lambda+1=0,\;\;\mu^{2}-c\mu+b=0
\]
respectively. Observe that
\begin{equation}
b<1\Rightarrow\lambda_{1}<\mu_{1}\label{mb}%
\end{equation}
Let us define functions $\varphi_{1}$ and $\varphi_{2}$ as follows:
\[
\varphi_{1}(t):=\left\{
\begin{array}
[c]{l}%
\frac{1}{2}e^{\lambda_{1}t},\;\;\;\;\;\;\quad t\leq0,\\
1-\frac{1}{2}e^{-\lambda_{1}t},\quad t>0
\end{array}
\right.  \;\;\;\;\varphi_{2}(t):=\left\{
\begin{array}
[c]{l}%
\frac{1}{2}e^{\mu_{1}t},\;\;\;\;\;\quad t\leq0,\\
1-\frac{1}{2}e^{-\mu_{1}t},\quad t>0
\end{array}
\right.
\]
observe that  $0<\varphi_{1}(t)<1\;$and similarly for
$0<\varphi_{2}(t)<1.$

\begin{lemma}
There exists a positive number $c^{\ast}=c^{\ast}(b,r)$ such that if
$c>c^{\ast}$, then $(\varphi_{1},\varphi_{2})^{T}$ defined as above
is a smooth monotone upper solution of the wave equation.
\end{lemma}
\begin{proof}
First, it is easily seen that
\begin{align*}
\varphi_{1}^{\prime}(t) &  =\left\{
\begin{array}
[c]{l}%
\frac{\lambda_{1}}{2}e^{\lambda_{1}t},\quad t\leq0,\\
\frac{\lambda_{1}}{2}e^{-\lambda_{1}t},\quad t>0
\end{array}
\right.  ,\quad\varphi_{1}^{\prime\prime}(t)=\left\{
\begin{array}
[c]{l}%
\frac{\lambda_{1}^{2}}{2}e^{\lambda_{1}t},\quad t\leq0,\\
\frac{-\lambda_{1}^{2}}{2}e^{-\lambda_{1}t},\quad t>0
\end{array}
\right.  \\
\varphi_{2}^{\prime}(t) &  =\left\{
\begin{array}
[c]{l}%
\frac{\mu_{1}}{2}e^{\mu_{1}t},\quad t\leq0,\\
\frac{\mu_{1}}{2}e^{-\mu_{1}t},\quad t>0
\end{array}
\right.  ,\quad\varphi_{2}^{\prime\prime}(t)=\left\{
\begin{array}
[c]{l}%
\frac{\mu_{1}^{2}}{2}e^{\mu_{1}t},\quad t\leq0,\\
\frac{-\mu_{1}^{2}}{2}e^{-\mu_{1}t},\quad t>0
\end{array}
\right.
\end{align*}
The second derivative exists almost everywhere. substituting the
above expressions into the first equation of (\ref{bz-we}) we have
for $t\leq0$
\begin{align*}
&  \left[
\varphi_{1}^{\prime\prime}(t)-c\varphi_{1}^{\prime}(t)+\varphi
_{1}(t)\right]  -r\varphi_{1}(t)(1-\varphi_{2}(t-c\tau))-\varphi_{1}^{2}(t)\\
&  =-re^{\lambda_{1}t}\left(  1-\varphi_{2}(t-c\tau)\right)  -\frac{1}%
{4}e^{2\lambda_{1}t}\leq 0
\end{align*}

since $\left(  1-\varphi_{2}(t-c\tau)\right)  >0.$

For the second component we have for $t\leq0$
\begin{align*}
&  \varphi_{2}^{\prime\prime}(t)-c\varphi_{2}^{\prime}(t)+b\varphi
_{1}(t)\left(  1-\varphi_{2}(t)\right)
=\frac{\mu_{1}^{2}}{2}e^{\mu_{1}t}-c\frac{\mu_{1}}{2}e^{\mu_{1}t}+\frac
{b}{2}e^{\lambda_{1}t}\left(  1-\frac{1}{2}e^{\mu_{1}t}\right)  \\
&  =\frac{b}{2}\left(  e^{\lambda_{1}t}-e^{\mu_{1}t}\right)  -\frac{b}%
{4}e^{\left(  \lambda_{1}+\mu_{1}\right)  t}%
\end{align*}
Thus by (\ref{mb})
\[
\lambda_{1}\leq\mu_{1}\Rightarrow\varphi_{2}^{\prime\prime}(t)-c\varphi
_{2}^{\prime}(t)+b\varphi_{1}(t)\left( 1-\varphi_{2}(t)\right)
\leq-\frac {b}{4}<0
\]
On the other hand if $t\geq0,$ then we have%
\begin{align*}
&  \left[
\varphi_{1}^{\prime\prime}(t)-c\varphi_{1}^{\prime}(t)+\varphi
_{1}(t)\right]  -r\varphi_{1}(t)(1-\varphi_{2}(t-c\tau))-\varphi_{1}^{2}(t)\\
&  =\left[  \frac{-\lambda_{1}^{2}}{2}e^{-\lambda_{1}t}-c\frac{\lambda_{1}}%
{2}e^{-\lambda_{1}t}+1-\frac{1}{2}e^{-\lambda_{1}t}\right]  -r\left(
1-\frac{1}{2}e^{-\lambda_{1}t}\right) \left(
1-\varphi_{2}(t-c\tau)\right)
-\varphi_{1}^{2}(t)\\
&  =-r\left(  1-\frac{1}{2}e^{-\lambda_{1}t}\right)  \left(
1-\varphi
_{2}(t-c\tau)\right)  -\lambda_{1}^{2}e^{-\lambda_{1}t}-\frac{1}%
{4}e^{-2\lambda_{1}t}\\
&  <0 .
\end{align*}
Similarly the second equation yields%
\begin{align*}
&  \varphi_{2}^{\prime\prime}(t)-c\varphi_{2}^{\prime}(t)+b\varphi
_{1}(t)\left(  1-\varphi_{2}(t)\right)  \\
&  =\frac{-\mu_{1}^{2}}{2}e^{-\mu_{1}t}-c\frac{\mu_{1}}{2}e^{-\mu_{1}%
t}+b\left(  1-\frac{1}{2}e^{-\lambda_{1}t}\right)
-b\varphi_{1}(t)\varphi
_{2}(t)\\
&  =-ce^{-\mu_{1}t}+\frac{b}{2}\left(
e^{-\mu_{1}t}-e^{-\lambda_{1}t}\right)
+\frac{b}{2}e^{-\mu_{1}t}+\frac{b}{2}e^{-\lambda_{1}t}-\frac{b}{4}e^{-\left(
\lambda_{1}+\mu_{1}\right)  t}\\
&  =\left[  \left(  b-c\right)
-\frac{b}{4}e^{-\lambda_{1}t}\right] e^{-\mu_{1}t},
\end{align*}
so, when $t\geq0$,
\[
b-c\leq0\Rightarrow\varphi_{2}^{\prime\prime}(t)-c\varphi_{2}^{\prime
}(t)+b\varphi_{1}(t)\left( 1-\varphi_{2}(t)\right)  <0
\]
To summarize the situation we have
\[
\left.
\begin{array}
[c]{c}%
b\leq c\\
\lambda_{1}\leq\mu_{1}%
\end{array}
\right\}
\Rightarrow\mathbb{\;\;}\varphi_{2}^{\prime\prime}(t)-c\varphi
_{2}^{\prime}(t)+b\varphi_{1}(t)\left( 1-\varphi_{2}(t)\right)
<0\;\;\;\;\forall t\in\mathbb{R,}%
\]
and by (\ref{mb}) the condition reduces to
\[
\left\{
\begin{array}
[c]{c}%
b\leq c\\
b<1
\end{array}
\right.  \Rightarrow\mathbb{\;\left\{
\begin{array}
[c]{c}%
b\leq c\\
\lambda_{1}\leq\mu_{1}%
\end{array}
\right.  \;}%
\]
which is satisfied if
\[
b<\min\left\{  1,c\right\}  .
\]
\end{proof}
Finally, since the procedure of constructing quasi-lower solutions
is similar for the predator-prey models we leave it to the reader.
To conclude this section we have
\begin{theorem}\label{the 3}
There exists a positive number $c^{\ast}=c^{\ast}(b,r)$ such that if
$c>c^{\ast}$, then Eq. (\ref{pde-bz}) has a wave front solution
$u(x,t)=\varphi_1(x+ct)$, $v(x,t)=\varphi_2(x+ct)$.
\end{theorem}


\begin{thebibliography}{99}

\bibitem{chilat} {C. Chicone, Y. Latushkin}, "Evolution
Semigroups in  Dynamical  Systems and Differential Equations",
Mathematical Surveys and Monographs, vol. 70, Providence, RO:
American Mathematical Society, 1999.

\bibitem{chohal} Shui Nee Chow, Jack K. Hale,  "Methods of Bifurcation
Theory". Grundlehren der Mathematischen Wissenschaften, {\bf 251}.
Springer-Verlag, New York-Berlin, 1982.
\bibitem{congar}  C. Conley and R. Gardner, An application of the
generalized Morse index to traveling wave solutions of a competitive
reaction-diffusion model, {\it Indiana Univ. Math. J.} {\bf 44}
(1984), 319-343.
\bibitem{dun1} S. R. Dunbar, Traveling wave solutions of diffusive Lotka-Volterra
equations, {\it J. Math. Biol.}, {\bf 17} (1983), 11-32.

\bibitem{dun2} S. R. Dunbar, Traveling wave solutions of diffusive Lotka-Volterra
equations: A hetero- clinic connection in $R^4$, {\it Trans. Amer.
Math. Soc.} {\bf 268} (1984), 557-594.

\bibitem{farhuawu} T. Faria, W. Huang, J. Wu, Travelling waves for delayed reaction – diffusion
equations with global response. {\it Proc. R. Soc. A}, {\bf 462}
(2006), 229 –261.

\bibitem{fis}  R. A. Fisher, The wave of advance of advantageous genes, {\it Ann.
Eugenics}, {\bf 7} (1937), 355-369.

\bibitem{gar}  R. Gardner, Existence of traveling wave solutions of predator-prey
systems via the connection index. {\it SIAM J. Appl. Math.} {\bf 44}
(1984), 56-79.
\bibitem{gar2}  R. Gardner, Existence and stability of
traveling wave solutions of competition models: A degree theoretic
approach, J. Differential Equations 44 (1982), 343-364.

\bibitem{gar3}R. Gardner, Review on Traveling Wave Solutions of Parabolic
Systems by A. I. Volpert, V. A. Volpert and V. A. Volpert. {\it
Bull. Am. Math. Soc.}, {\bf 32} (1995), 446-452.

\bibitem{hal}J. Hale, Theory of Functional Differential
Equations,  Acad. Press, New York, 1977.


\bibitem{hen}
D. Henry, {"Geometric Theory of Semilinear Parabolic Equations"},
Lecture Notes in Math., Springer-Verlag, Berlin-New York, 1981.



\bibitem{hinnaiminshi}
Y. Hino, T. Naito, N.V. Minh, J.S. Shin, "Almost Periodic
Solutions of Differential Equations in Banach Spaces". Taylor \&
Francis, London - New York, 2002.

\bibitem{hualurua} J. Huang, G. Lu, S. Ruan,
Existence of traveling wave solutions in a diffusive predator-prey
model. {\it J. Math. Biol.} {\bf 46} (2003), 132--152.

\bibitem{huazou} Jianhua Huang, Xingfu Zou,  Existence of
traveling wavefronts of delayed reaction diffusion systems without
monotonicity. {\it Discrete Contin. Dyn. Syst.} {\bf 9} (2003),
925--936.


\bibitem{kolpetpis}  A. N. Kolomgorov, I. G. Petrovskii, and N. S. Piskunov, Study
of a diffusion equation that is related to the growth of a quality
of matter, and its application to a biological problem, {\it Byul.
Mosk. Gos. Univ. Ser. A Mat. Mekh.} {\bf 1} (1937), 1-26.



\bibitem{kan}  Ya. I. Kanel, Existence of a
traveling-wave solution of the Belousov-Zhabotinskii system, {\it
Differentsial'nye Uravneniya}, {\bf 26} (1990), 652-660.

\bibitem{kap}  A. Ya. Kapel, Existence of traveling-wave type solutions for
the Belousov-Zhabotinskii system equations, {\it Sibirsk. Mat.
Zh.}, {\bf 32} (1991), 47-59.

\bibitem{liunguminvu} J. Liu, G. Nguerekata, Nguyen Van Minh, Vu Quoc Phong, Bounded solutions of parabolic equations in continuous
function spaces. To appear in {\it Funkciolaj Ekvacioj}.


\bibitem{ma} Shiwang Ma, Traveling Wavefronts for Delayed Reaction-Diffusion
Systems via a Fixed Point Theorem. {\it Journal of Differential
Equations}, {\bf 171} (2001), 294-314.

 \bibitem{mazou} Shiwang Ma,  Xingfu Zou, Existence, uniqueness and stability of travelling waves in a discrete reaction-diffusion
 monostable equation with delay. {\it J. Differential Equations}, {\bf 217} (2005), 54--87.


\bibitem{mis} K. Mischaikow,  Travelling waves for a cooperative and a competitive-cooperative system.
in "{\it Viscous profiles and numerical
 methods for shock waves}" (Raleigh, NC, 1990), 125--141, SIAM, Philadelphia, PA, 1991.

\bibitem{mishut} K. Mischaikow,  V. Hutson, Travelling waves for
mutualist species. {\it SIAM J. Math. Anal.} {\bf 24} (1993),
987--1008.

\bibitem{misrei} K. Mischaikow, J. F. Reineck, James F.
Travelling waves in predator-prey systems. {\it  SIAM J. Math.
Anal.} {\bf 24} (1993), 1179--1214.


\bibitem{murnaimin}
S. Murakami, T. Naito, N. van Minh, {Evolution semigroups and sums
of commuting operators: A new approach to the admissibility theory
of function spaces}, {\it J. Differential Equations}, {\bf 164}
(2000), 240-285.

\bibitem{naimin}T. Naito and N. van Minh,  Evolutions semigroups and spectral
criteria for almost periodic solutions  of periodic evolution
equations, {\it Journal of Differential Equations}, {\bf 152}
(1999), 358-376.

\bibitem{mur}J. D. Murray, "Mathematical Biology". Berlin Heidelberg New York:
Springer, (1993).

\bibitem{sch}  K. Schaaf, Asymptotic behavior and traveling wave solutions for
parabolic functional differential equations, {\it Trans. Amer.
Math. Soc.}, {\bf 302} (1987), 587-615.

\bibitem{tro}  W. C. Troy, The existence of traveling wavefront solutions of a
model of the Belousov Zhabotinskii reaction, {\it J. Differential.
Equations}, {\bf 36} (1980), 89-98.

\bibitem{3vol}  A. I. Volpert, V. A. Volpert, and
V. A. Volpert, "Traveling Wave Solutions of Parabolic Systems",
Translations of Mathematical Monographs, Vol. 140, Amer. Math.
Soc., Providence, 1994.

\bibitem{wanlirua} Zhi-Cheng Wang, Wan-Tong Li,  Shigui Ruan, Travelling wave
fronts in reaction-diffusion systems with spatio-temporal delays.
{\it J. Differential Equations}, {\bf 222} (2006), 185--232.

\bibitem{wuzou} J. Wu and X. Zou, Traveling wave fronts of
reaction-diffusion systems with delay, {\it Journal of Dynamics and
Differential Equations}, {\bf 13} (2001), 651-687.


\bibitem{tanfif}  M. M.
Tang and P. Fife, Propagating fronts for competing species equations
with diffusion, {\it Arch. Rational Mech. Anal.} {\bf 73} (1980),
69-77.

\bibitem{wu}J. Wu, "Theory and Applications of Partial Functional Differential
Equations", Applied Math. Sci. 119, Springer, Berlin-New York,
1996.

\bibitem{yewan}  Q. Ye and M. Wang, Traveling wavefront solutions of Noyes-field
system for Belousov Zhabotinskii reaction, {\it Nonlinear Anal.},
{\bf 11} (1987), 1289-1302.

\bibitem{zouwu}  X. Zou and J. Wu, Existence of traveling
wavefronts in delayed reaction-diffusion system via monotone
iteration method, {\it Proc. Amer. Math. Soc.}, {\bf 125} (1997),
2589-2598.

\bibitem{zha}  Xiao-Qiang Zhao, "Dynamical systems in population biology". CMS
Books in Mathematics/Ouvrages de Mathématiques de la SMC, 16.
Springer-Verlag, New York, 2003.

\end{thebibliography}
\end{document}